\documentstyle{amsart}

\newcommand{\cA}{{\mathcal A}}
\newcommand{\cS}{{\mathcal S}}
\newcommand{\cO}{{\mathcal O}}
\newcommand{\cL}{{\mathcal L}}
\newcommand{\cB}{{\mathcal B}}

\newcommand{\cG}{{\mathcal G}}
\newcommand{\cV}{{\mathcal V}}
\newcommand{\cH}{{\mathcal H}}
\newcommand{\cN}{{\mathcal N}}
\newcommand{\fg}{{\mathfrak g}}
\newcommand{\fL}{{\mathfrak L}}

\newcommand{\fl}{{\mathfrak l}}
\newcommand{\fb}{{\mathfrak b}}
\newcommand{\fn}{{\mathfrak n}}
\newcommand{\uc}{\underline{c}}
\newcommand{\codim}{\mbox{codim}}
\newcommand{\Ad}{\mbox{Ad}}
\newcommand{\Lie}{\mbox{Lie}}
\newcommand{\Hom}{\mbox{Hom}}
\newcommand{\Ext}{\mbox{Ext}}

\newcommand{\BZ}{{\mathbb Z}}
\newcommand{\BN}{\Bbb{N}}
\newcommand{\BC}{{\mathbb C}}

\newcommand{\BQ}{{\mathbb Q}}
\newcommand{\bX}{{\bf X}}

\newcommand{\oo}{\overline{\phantom{A}}}

\newcommand{\Lm}{{\bf Lemma.}}
\newcommand{\Pf}{{\bf Proof. }}
\newcommand{\Rem}{{\bf Remark.}}
\newcommand{\Co}{{\bf Corollary.}}
\renewcommand{\l}{\langle}
\newcommand{\r}{\rangle}

\title{On the equivariant $K-$theory of the nilpotent cone}
\address{Independent Moscow University, 11 Bolshoj Vlasjevskij per.,
 Moscow
121002 Russia}
\date{October, 1999}
\email{ostrik@@mccme.ru}
\author{Viktor Ostrik}

\begin{document}
\maketitle
\section{Introduction}
Let $G$ be a simple simply connected algebraic group over the complex numbers.
Let $\fg$ be the Lie algebra of $G$. Let $\cN \subset \fg$ be the nilpotent cone.
Let $\cG =G\times \BC^*$. Let us consider the action of $\cG$ on $\fg$ given by
the rule $(g,z)x=z^{-2}\Ad (g)x$. The nilpotent cone $\cN$ is invariant
under this action. The aim of this note is to make some conjectures on the 
equivariant
$K-$group $K_{\cG}(\cN)$, see e.g. \cite{CG}. Namely, we introduce a
"Kazhdan-Lusztig type" canonical basis of $K_{\cG}(\cN)$ over the 
representation
ring of $\BC^*$, parametrized by dominant weights of $G$. We conjecture that
this basis is close to the basis consisting of irreducible $G-$equivariant
bundles on nilpotent orbits. This would give us a bijection between two sets:
\{ dominant weights for $G$\} and \{ pairs consisting of a nilpotent orbit
$\cO$ and an irreducible $G-$equivariant bundle on $\cO$\}. Such a bijection
appeared in the work of G.Lusztig on the asymptotic affine Hecke algebra, see
\cite{L4} IV 10.8. We conjecture that our (conjectural) bijection coincides
with Lusztig's. We also conjecture that some specific elements of $K_{\cG}(\cN)$
closely related with irreducible local systems on nilpotent orbits belong to our
basis. All our Conjectures are motivated by the study of
the cohomology of {\em quantized tilting modules}. So this note should be
considered as a generalization of Humphreys' Conjecture \cite{Hu}.

This note is based on the idea of George Lusztig. I learnt
this idea from Michael Finkelberg.
I am deeply grateful to Roman Bezrukavnikov
for extremely useful conversations while preparing this note.
I wish to thank David Vogan who pointed out a gap in the first version of this
note. Thanks also due to Jim Humphreys for valuable suggestions. Finally
I would like to acknowledge Harvard University for its hospitality while this
note was being written. 

\section{Canonical basis of $K_{\cG}(\cN)$}
\subsection{} \label{def}
Let $B\subset G$ be a Borel subgroup and let $\cB=G/B$ be the flag variety.
Let $\fb =\Lie (B)\subset \fg$ and let $\fn \subset \fg$ be the nilpotent
radical of $\fb$. It is well known that the cotangent bundle $T^*\cB$ is naturally
isomorphic to $G\times_B\fn$ and the map ({\em Springer resolution}) $s:
G\times_B\fn \to \cN,\; (g,n)\mapsto \Ad (g)n$ is a resolution of singularities
of $\cN$, see e.g. \cite{CG}.
This map is $\cG-$equivariant with respect to the $\cG-$action on $G\times_B\fn$
given by $(g,z)(g_1,n)=(gg_1,z^{-2}\Ad (g)n)$.

Let $X=\Hom (B,\BC^*)$ be the weight lattice of $G$ and let $X_+$ be the set of
dominant weights. For any $\lambda \in X_+$ let $V_{\lambda}$ denote the
irreducible representation of $G$ with highest weight $\lambda$. We will
also consider $V_{\lambda}$ as a $\cG-$module via projection $pr_1:G\times
\BC^*\to G$.

For any $\lambda \in X$ one associates the line bundle $\cL_{\lambda}$
on $\cB$ (see \cite{CG}). Let $\pi :G\times_B\fn \to \cB$ be the natural projection.
For any $\lambda \in X_+$ and $i>0$ we have $R^is_*\pi^*\cL_{\lambda}=0$
(Andersen-Jantzen vanishing), see
\cite{AJ,Br}; in this case we will call the sheaf $\widetilde{AJ}(\lambda)=
s_*\pi^*\cL_{\lambda}$ an {\em Andersen-Jantzen sheaf} (or AJ-{\em sheaf}).
Any AJ-sheaf $\widetilde{AJ}(\lambda)$ is endowed with the natural structure
of equivariant $\cG-$sheaf.

Let $R(\cG)$ be the representation ring of $\cG$. The ring $R(\cG)$ acts on
$K_{\cG}(V)$ for any $\cG-$variety $V$. Let $v\in R(\cG)$ correspond to the
one-dimensional representation of $\cG$ given by the second projection
$pr_2:\cG =G\times \BC^*\to \BC^*$. Let $\cA =\BZ [v,v^{-1}]\subset R(\cG)$ be the
subring of $R(\cG)$ generated by $v$ and $v^{-1}$.

Let $\tilde R(\cG)$ denote the $\cA-$module of formal linear combinations $\sum_{\lambda \in X_+}
k_{\lambda}V_{\lambda}, k_{\lambda}\in \cA$ (in general an infinite sum).
The space $\tilde R(\cG)$ is endowed with an obvious $R(\cG)-$structure.
For any $\cG-$equivariant sheaf $\cS$ on $\cN$, its global sections $\Gamma (\cS)$
form a $\cG-$module in a natural way and $\Hom_G(V_{\lambda}, \Gamma (\cS))$ is
finite dimensional for any $\lambda \in X_+$. Since $\cN$ is affine the functor
$\Gamma$ is exact and we obtain a well-defined map $\Gamma :K_{\cG}(\cN)\to \tilde R(\cG)$.
Clearly this map is $R(\cG)-$linear. Note that for any $\lambda \in X_+$
$$
\Gamma (\widetilde{AJ}(\lambda))\in V_{\lambda}+\sum_{\lambda'\in X_+}v\BZ [v]V_{\lambda'}. \leqno (a)
$$

\subsection{} \label{basis}
\Lm (R.Bezrukavnikov) {\em The classes $[\widetilde{AJ}(\lambda)]\in
K_{\cG}(\cN), \lambda \in X_+$ form an $\cA -$basis of $K_{\cG}(\cN)$.}

\Pf Let us prove that $[\widetilde{AJ}(\lambda)]$ are linearly independent.
Suppose that \linebreak $\sum_{\lambda \in X_+}a_{\lambda}[\widetilde{AJ}(\lambda)]=0,\;
a_{\lambda}\in \cA$. We may assume that $a_{\lambda}\in \BZ [v]$ and
$a_{\lambda}\not \in v\BZ [v]$ for at least one $\lambda$. Applying $\Gamma$ to the
both sides of this equality and using \ref{def} (a) we obtain a contradiction.

Let us prove that the map $s_*: K_{\cG}(T^*\cB)\to K_{\cG}(\cN)$ is surjective.
To this end let us show that for any $\cG-$equivariant sheaf $\cS$ on $\cN$
there exists $\alpha \in K_{\cG}(T^*\cB)$ such that $[\cS]-s_*\alpha$ has a
strictly smaller support than $[\cS]$ (then we are done by devissage).
We can assume that the support of $[\cS]$ is the closure $\bar \cO$ of a
nilpotent orbit $\cO$. But it is well known (see e.g. \cite{Mc1}) that
$\bar \cO$ admits a resolution of singularities $r: X\to \bar \cO$ where
$X$ is some $G-$equivariant subbundle of the cotangent bundle of some partial 
flag variety. Clearly the
support of $[\cS]-r_*[r^*\cS]$ is contained in $\bar \cO -\cO$. Finally, one
shows using the Koszul complex that the image of $K_{\cG}(X)$ under $r_*$ is
contained in the image of $K_{\cG}(T^*\cB)$ under $s_*$. The surjectivity is
proved. 

As it is well known the sheaves $\pi^*\cL_{\lambda}, \lambda \in X$ form
an $\cA -$basis in $K_{\cG}(T^*\cB)$, see e.g. \cite{CG}. Now let
$\lambda \in X\setminus X_+$. Let $\alpha_i$ be a simple root such that
$\l \lambda, \alpha_i^{\vee}\r <0$. A simple $SL_2-$calculation (see e.g.
\cite{Br1} 3.15) shows that
$$
[s_*\pi^*\cL_{\lambda}]=\left\{ \begin{array}{cc}v^2[s_*\pi^*\cL_{s_{\alpha_i}\lambda}]&\mbox{if}\;
\l \lambda, \alpha_i^{\vee}\r =-1;\\
-[s_*\pi^*\cL_{s_{\alpha_i}\lambda-\alpha_i}]+
v^2[s_*\pi^*\cL_{s_{\alpha_i}\lambda}]+v^2[s_*\pi^*\cL_{\lambda +\alpha_i}]&\mbox{if}\;
\l \lambda, \alpha_i^{\vee}\r \le -2.\end{array}\right.
$$
The Lemma follows.

\subsection{}
Let $W$ be the Weyl group of $G$ and let $\nu$ be the number of positive roots in $G$.
Let $l:W\to \BN$ be the length function and let $w_0\in W$ be the longest element (then $\nu =l(w_0)$).
For any $\lambda \in X$ let $W_{\lambda}\subset W$ be the stabilizer of $\lambda$
in $W$ and let $\nu_{\lambda}$ be the length of the longest element $w_{\lambda}$
of $W_{\lambda}$.

It follows from the Lemma \ref{basis} that the classes $AJ(\lambda)=(-v)^{\nu -\nu_{\lambda}}
[\widetilde{AJ}(\lambda)], \lambda \in X_+$ form a $\BZ [v,v^{-1}]-$basis of
$K_{\cG}(\cN)$. We will call $\{ AJ(\lambda)\}$ the {\em Andersen-Jantzen basis}.

\subsection{} \label{sphere}
Let $Z=T^*\cB\times_{\cN}T^*\cB$ be the Steinberg variety, see \cite{L6}.
Let $\cH$ be the affine Hecke algebra (over $\BZ [v,v^{-1}]$) associated with $G$, see \cite{L6}.
We identify $K_{\cG}(Z)$ and $\cH$ via the isomorphism constructed in \cite{L6} (see also \cite{CG}).
Let $\hat W^a$ be the (extended) affine Weyl group as defined in \cite{L6} 1.7.
The Weyl group $W$ is embedded in $\hat W^a$ and the set of double cosets
$W\setminus \hat W^a/W$ is identified with the set $X_+$, see \cite{L4}; for any
$w\in \hat W^a$ let $\lambda_w\in X_+$ denote its double coset. Let $l(\lambda_w)$
denote the length of the shortest element in $\lambda_w$.
Let $\tilde T_w, w\in \hat W^a$ be the basis of $\cH$ as defined in \cite{L6} 1.8.
The natural projection $st:Z\to \cN$ induces homomorphism $st_*:K_{\cG}(Z)\to
K_{\cG}(\cN)$.

\Lm {\em We have an equality 
$st_*(\tilde T_w)=(-v)^{-l(w)+l(\lambda_w)}AJ(\lambda_w)$.
In particular, the map $st_*$ is surjective.}

\Pf The map $st$ can be factorized in two ways: $st=s\cdot pr_1$ and
$st=s\cdot pr_2$ where $pr_i, i=1,2$ are two projections $Z\to T^*\cB$. It
follows (see \cite{L6} 7.25 and 7.19) that $st_*(T_{s_iw})=-v^{-1}st_*(T_w)$ and
$st_*(T_{ws_i})=-v^{-1}st_*(T_w)$.
Finally, we note that our formula follows from definitions in \cite{L6} for any
translation by dominant weight considered as element of $\hat W^a$.

{\bf Remark.} Let $Cent(\cH)$ be the center of $\cH$ (see e.g. \cite{CG} for
its description). It {\em is not} true that the map $st_*: Cent(\cH)\to K_{\cG}(\cN)$
is surjective, in fact its image consists of trivial bundles on $\cN$ with
possibly nontrivial $\cG-$structure. But it is a consequence of Lemmas
\ref{basis} and \ref{sphere} that the map $st_*: Cent(\cH)\otimes_{\cA}\BQ (v)\to
K_{\cG}(\cN)\otimes_{\cA}\BQ (v)$ is surjective.

\subsection{}
Let $\oo :\cH \to \cH$ be the Kazhdan-Lusztig involution of the ring $\cH$, see \cite{L6} 1.8.
It follows from Lemma \ref{sphere} that the kernel of the map $st_*$ is invariant
under $\oo $. So, we obtain a well defined involution $\oo :K_{\cG}(\cN)\to
K_{\cG}(\cN)$. We call this map Kazhdan-Lusztig involution.

Let $\omega :\fg \to \fg$ be an opposition (see \cite{L6} 9.8).
Let $D: K_{\cG}(\cN)\to K_{\cG}(\cN)$ be the Serre-Grothendieck duality map,
see \cite{L6} 6.10.

\Lm {\em The Kazhdan-Lusztig involution $\oo :K_{\cG}(\cN)\to K_{\cG}(\cN)$
equals $v^{-2\nu}D\omega^*$.}

\Pf This immediately follows from \cite{L6} 9.12 if we note that restriction of
$st_*$ to the center of the Hecke algebra is surjective after tensoring with $\BQ (v)$.

\subsection{}
We say that $x\in K_{\cG}(\cN)$ is selfdual if $\bar x=x$.

\Lm {\em For any $\lambda \in X_+$ there exists a unique selfdual element
$C(\lambda)\in K_{\cG}(\cN)$ such that $C(\lambda)\in AJ(\lambda)+
\sum_{\mu \in X_+}v^{-1}\BZ [v^{-1}]AJ(\mu)$. The elements $C(\lambda)$ form a
basis of $K_{\cG}(\cN)$.}

\Pf Let $c_w'\in \tilde T_w+\sum_{w'<w}v^{-1}\BZ [v^{-1}]\tilde T_{w'}, \bar c_w'=c_w',
w\in \hat W^a$ be the Kazhdan-Lusztig basis of $\cH$, see \cite{L6} 1.5, 1.8.
For any $\lambda \in X_+=W\setminus \hat W^a/W$ let $m_{\lambda}\in \lambda$ be
the shortest element.
We set $C(\lambda)=st_*(c_{m_{\lambda}}')$. The unicity follows from the
existence in a standard way, see e.g. \cite{So} 2.4.

\Rem Let $C(\lambda')=\sum_{\lambda}b_{\lambda ,\lambda'}AJ(\lambda)$ where
$b_{\lambda ,\lambda'}\in \BZ [v^{-1}]$.
The polynomials $b_{\lambda ,\lambda'}$ appeared
in the work of G.Lusztig \cite{L2}.
The idea that the matrix $b_{\lambda, \lambda'}$ or rather
its inverse should have representation theoretic meaning is due to
Ivan Mirkovi\'c. We believe that this note is a step in this direction.

\Co {\em We have $st_*(c_w')=0$ unless $w=m_{\lambda_w}$ in which
case $st_*(c_w')=C(\lambda_w)$.}

\section{Conjectures}
In this section we formulate a number of conjectures on the basis $\{ 
C(\lambda)\}$.

\subsection{}
For any $C\in K_{\cG}(\cN)$ one defines its {\em support} as the complement to
the union of all open $j:U\hookrightarrow \cN$ such that $j^*C=0$. Clearly, the
support of any $C\in K_{\cG}(\cN)$ is a closed $G-$invariant subset of $\cN$.

{\bf Conjecture 1.} {\em The support of any $C(\lambda)$ is irreducible.}

In other words the support of any element $C(\lambda)$ is the closure of some
nilpotent orbit $\cO =\cO_{\lambda}$. Let $j_{\cO}:\cO \hookrightarrow \cN$
be the inclusion. The element $j_{\cO_{\lambda}}^*C(\lambda)\in
K_{\cG}(\cO_{\lambda})$ is well defined.

{\bf Conjecture 2.} {\em The class $\pm j_{\cO_{\lambda}}^*C(\lambda)$ is
represented by an irreducible $\cG-$equivariant bundle $\cV_{\lambda}$ on $\cO_{\lambda}$.}

Let us choose a set $\{ e_{\cO}\}$ of representatives of all nilpotent orbits.
Let $C_{\cG}(e_{\cO})$ (resp. $C_{G}(e_{\cO})$) be the centralizer of $e_{\cO}$
in $\cG$ (resp. in $G$). Irreducible representations
of $C_{\cG}(e_{\cO})$ (resp. $C_{G}(e_{\cO})$) factor through the quotient
$C^{red}_{\cG}(e_{\cO})$ (resp. $C^{red}_{G}(e_{\cO})$) of $C_{\cG}(e_{\cO})$
(resp. $C_{G}(e_{\cO})$) by its unipotent radical. The exact sequence
$$
1\to C_{G}(e_{\cO})\to C_{\cG}(e_{\cO})\to \BC^*\to 1 \eqno (1)
$$
where the first map is obvious inclusion and the second one is restriction of
second projection $pr_2:\cG \to \BC^*$ induces an exact sequence
$$
1\to C^{red}_{G}(e_{\cO})\to C^{red}_{\cG}(e_{\cO})\to \BC^*\to 1. \eqno (2)
$$
We remark that the last sequence canonically splits. Namely, let $\phi_{\cO}: SL_2(\BC)\to
G$ be a homomorphism such that $d\phi_{\cO}\left( \begin{array}{cc}0&1\\ 
0&0\end{array}
\right) =e_{\cO}$, which exists by Jacobson---Morozov Theorem. Then
$\psi_{\cO}: \BC^*\to \cG,\; \psi_{\cO}(z)=(\phi_{\cO}
\left( \begin{array}{cc}z^{-1}&0\\ 0&z\end{array}\right) , z)$ is a desired
splitting. It is not canonical for the sequence (1) since it requires choice
of $\phi_{\cO}$ but all choices give the same splitting of the sequence (2).

Taking the stalk at $e_{\cO}$ defines an equivalence
of categories \{ $\cG-$equivariant bundles on $\cO$\} and \{ representations
of $C_{\cG}(e_{\cO})$\} , see e.g. \cite{CG}. So the Conjecture 2 gives us
an irreducible representation $\rho_{\lambda}'$ of $C_{\cG}(e_{\cO_{\lambda}})$
attached to the bundle $\cV_{\lambda}$. By the above we can consider $\rho_{\lambda}'$
as representation of $C^{red}_{\cG}(e_{\cO_{\lambda}})$. Let $a(\cO )=\frac{1}{2}\codim_{\cN}\cO$.
We expect that the group $\psi_{\lambda}(\BC^*)$ acts on $\rho_{\lambda}'$ by
dilatations $z\mapsto z^{-a(\cO_{\lambda})}Id$. So the representation
$\rho_{\lambda}'$ is completely defined by its restriction $\rho_{\lambda}$
to $C^{red}_G(e_{\cO_{\lambda}})$.

The group $C^{red}_G(e_{\cO_{\lambda}})$ contains a canonical central involution
$\phi_{\cO_{\lambda}}\left( \begin{array}{cc}-1&0\\ 0&-1\end{array}\right)$. It
acts on the irreducible representation $\rho_{\lambda}$ by $\pm Id$. We expect
that the
sign here is the same as the sign in the Conjecture 2.

Conjectures 1 and 2 provide a map $\fL: \lambda \mapsto (\cO_{\lambda},
\rho_{\lambda})$ from the set of dominant weights to the set of pairs
\{ nilpotent orbit $\cO$, irreducible representation of $C^{red}_G(e_{\cO})$\} .

{\bf Conjecture 3.} {\em The map $\fL$ is a bijection. Moreover, this bijection
coincides with Lusztig's defined in \cite{L4} IV, 10.8.}

In particular the orbit $\cO_{\lambda}$ can be determined as follows. We take
an element $m_{\lambda}\in \hat W^a$ and find the {\em two-sided cell}
$\uc =\uc_{\lambda}\subset \hat W^a$
containing $m_{\lambda}$, see \cite{L4}. The main result of \cite{L4} IV is a
bijection $\fl$ between the set of two sided cells and the set of nilpotent orbits.
Humphreys' Conjecture \cite{Hu} predicts that $\cO_{\lambda}=\fl (\uc_{\lambda})$.

Now let $\cO$ be a nilpotent orbit and let $\BC_{\cO}$ be the
trivial one dimensional bundle on $\cO$. The sheaf $\BC_{\cO}$ has an obvious
$\cG-$equivariant structure. Let $j:\cO \to \cN$ be the inclusion. Consider
the sheaf $j_*\BC_{\cO}$.

{\bf Conjecture 4.} {\em The element $v^{-a(\cO)}[j_*\BC_{\cO}]\in K_{\cG}(\cN)$
is of the form $C(\lambda)$ for some $\lambda =\lambda_{\cO}$. Moreover,
the element $d_{\cO}=m_{\lambda_{\cO}}$ is a distinguished involution in
$\hat W^a$, see \cite{L4} II.}

Remark that self-duality of the element $v^{-a(\cO)}[j_*\BC_{\cO}]$ follows
from the Theorem of Hinich and Panyushev \cite{Hi,P} stating that normalization
of the closure of any nilpotent orbit has rational singularities. Furthermore, let
$\rho$ be an irreducible representation of some $C_{G}(e_{\cO})$ which
factors through finite quotient of $C_{G}(e_{\cO})$. Using splitting of the
exact sequence (2) we extend $\rho$ to the representation of $C_{\cG}(e_{\cO})$
trivial on $\psi_{\cO}(\BC^*)$. Let $\cV_{\rho}$ be the corresponding
$\cG-$equivariant bundle on $\cO$. We expect that the class
$v^{-a(\cO)}[j_*\cV_{\rho}]\in K_{\cG}(\cN)$ also is of the form $C(\lambda)$
for some $\lambda =\lambda_{\rho}$. The self-duality of this element should be
a consequence of \cite{BrD} 6.3. We note that by no means is such a statement
true for general $\rho$, see \ref{sl3} below.

\subsection{Tilting modules} This subsection is devoted to the explanation
of the connection of our Conjectures with the theory of tilting modules.
In fact the Conjectures above were
motivated by the study of cohomology of tilting modules for quantum groups.
Tilting modules provide a lifting of these Conjectures from the K-theoretical
level to the level of categories. We refer the reader
to \cite{An} for the definition and basic properties of tilting modules.

Let $U$ be the quantum group over $\BC$ with the same root datum as the group
 $G$ at a primitive $l-$th
root of unity where $l$ is an odd number (prime to 3 if $G$ is of type $G_2$)
greater than Coxeter number of $G$. Let $u\subset U$ be the small quantum group.
Let $\BC$ denote the trivial representation of $U$. The Ginzburg--Kumar Theorem
\cite{GK} states that odd cohomology $H^{odd}(u, \BC)$ vanishes and the
graded algebra of even cohomology $H^{ev}(u, \BC)$ is isomorphic to the algebra
$\BC [\cN ]$ of regular functions on the nilpotent cone $\cN$ (grading on the
latter algebra comes from $\BC^*$-action). Moreover, the natural $G$-actions on
both algebras are the same under this isomorphism. For any finite dimensional
$U-$module $M$ the cohomology $H^{\bullet}(u, M)$ is a module over $H^{\bullet}(u, \BC )$
via cup-product. This module is finitely generated, see {\em loc. cit.} So we
 can
identify $H^{\bullet}(u, M)$ with a $\cG-$equivariant coherent sheaf on $\cN$, see \cite{GK}.
Further, we can attach to $M$ the class of its Euler characteristic $\chi (M)=
[H^{ev}(u, M)]-[H^{odd}(u, M)]\in K_{\cG}(\cN)$.

Now let $\bX$ be the weight lattice of $U$. Of course this lattice is isomorphic
to $X$ but we prefer to distinguish two lattices (in fact it is natural to
identify $X$ with the sublattice $l\bX \subset \bX$). Let us define a {\em dot-action}
$(w,x)\mapsto w\cdot x$ of the group $\hat W^a$ on $\bX$ as follows. Recall
that the group $\hat W^a$ is canonically isomorphic to the semidirect product
of $W$ with $X$, see e.g. \cite{L4} IV, 1.6. Now for any $w=\lambda v\in \hat W^a$
with $\lambda \in X,\; v\in W$ we set $w\cdot x=v(x+\varrho)-\varrho +l\lambda$.
Here $\varrho \in \bX$ is the half sum of positive roots.

For any dominant $x\in \bX$
let $T(x)$ denote an indecomposable tilting $U-$module with highest weight $x$
(it is unique up to nonunique isomorphism). In most cases the cohomology
$H^{\bullet}(u, T(x))$ vanishes. First of all the Linkage Principle shows
that $H^{\bullet}(u, T(x))\ne 0$ implies that $x\in \hat W_a\cdot 0$.
Further, we claim that $H^{\bullet}(u, T(w\cdot 0))\ne 0$ implies that
$w$ is minimal length element in its double coset $WwW\subset \hat W^a$.
Indeed, $w$ has minimal length in its coset $wW$ since $w\cdot 0$ is dominant.
If $w$ is not of minimal length in coset $Ww$ then there exists a simple
reflection $s\in W$ such that $l(sw)<l(w)$. Let $\Theta_s$ be the wall-crossing
functor corresponding to $s$, see e.g. \cite{So}. It easy to see that $T(w\cdot 0)$
is a direct summand of $\Theta_sT(sw\cdot 0)$. By adjointness of functors we
have:
$$
H^{\bullet}(u, \Theta_sT(sw\cdot 0))=\Ext_u^{\bullet}(\BC ,\Theta_sT(sw\cdot 0))=
\Ext_u^{\bullet}(\Theta_s\BC ,T(sw\cdot 0))=0
$$
since $\Theta_s\BC =0$. So the condition $H^{\bullet}(u, T(x))\ne 0$ implies
that $x=m_{\lambda}\cdot 0$ for some $\lambda \in X_+$. This condition is also
sufficient as shown by the following

{\bf Observation.} {\em Let $For: K_{\cG}(\cN)\to K_G(\cN)$ be the forgetting
map. Then $For (\chi (T(m_{\lambda}\cdot 0))=For (C(\lambda))$.}

This observation is a consequence of the (quantum version of) the main result 
of \cite{AJ}, Soergel's formula for characters of quantum tilting modules 
\cite{So}, and additivity of Euler characteristic. Indeed, Soergel's 
formula is in terms
of the $\cH-$module $\cH \otimes_{\cH_f} \cL (-v^{-1})$ where $\cH_f$ is Hecke
algebra of $W$ and $\cL (-v^{-1})$ is the one dimensional $\cH_f-$module
corresponding to the sign representation of $W$. Arguing as in the proof of
the Lemma \ref{sphere} we identify $\cH \otimes_{\cH_f} \cL (-v^{-1})$ with
$K_{\cG}(T^*\cB )$. We omit further details.

Conjectures 1---4 have their tilting versions. So we expect that the
sheaf $H^{\bullet}(u, T(m_{\lambda}\cdot 0))$ has irreducible support. In fact
we have two sheaves $H^{ev}(u, T(m_{\lambda}\cdot 0))$ and
$H^{odd}(u, T(m_{\lambda}\cdot 0))$ and we expect that their supports are related
by strict inclusion. The parity of the biggest cohomology sheaf
$H^{big}(u, T(m_{\lambda}))$ is determined by the sign in Conjecture 2 ($+$
corresponds to $H^{ev}$ and $-$ corresponds to $H^{odd}$). The sheaf
$H^{big}(u, T(m_{\lambda}))$ restricted to $\cO_{\lambda}$ should be equal to
$\cV_{\rho_{\lambda}}$. This picture is a generalization of Humphreys' Conjecture
on support varieties of tilting modules, see \cite{Hu}.

For any nilpotent orbit $\cO$ there exists a unique distinguished involution
$d_{\cO}$ such that $d_{\cO}$ is of minimal length in the double coset $Wd_{\cO}W$
and $d_{\cO}$ is contained in the cell $\uc$ with $\fl (\uc )=\cO$. The tilting
counterpart of Conjecture 4 states that the cohomology
$H^{\bullet}(u, T(d_{\cO}\cdot 0))$ vanishes in odd degrees, has a natural
structure of graded commutative algebra and is isomorphic as algebra to
the algebra $\BC [\cO ]$ of regular functions on $\cO$. This is true at least for
the regular nilpotent orbit (by Ginzburg-Kumar Theorem), for the trivial nilpotent
orbit (by \cite{An}) and the subregular nilpotent orbit (by \cite{O}). We also
expect that the parity vanishing holds in all cases when $\rho_{\lambda}$ is
trivial on the connected component of $C_G(e_{\cO_{\lambda}})$.

\section{Examples}
In this section we describe various cases when we were able to check some of
our Conjectures. It will be convenient to use the notations $e^{\lambda}=
[s_*\pi^*\cL_{\lambda}]$ for any $\lambda \in X$ (in particular if $\lambda \in X_+$
then $e^{\lambda}=[\widetilde{AJ}(\lambda)]$). The formula $\overline{e^{\lambda}}=
e^{w_0\lambda}$ which is a consequence of \cite{L6} 1.22 is very useful.

\subsection{$SL_2$} In this case $X=\BZ$ and $X_+=\BZ_{\ge 0}$. It is easy to
calculate that $C(0)=AJ(0)=e^0, C(1)=AJ(1)=e^1, C(2)=AJ(2)+v^{-1}AJ(0)=v^{-1}(e^0-v^2e^2),
C(n)=AJ(n)-v^{-2}AJ(n-2)=v^{-1}(e^{n-2}-v^2e^n)$ for $n\ge 3$. The support of
$C(0)$ and $C(1)$ is full nilpotent cone and the support of $C(n), n\ge 2$ is
the trivial nilpotent orbit. The element $C(0)$ (resp. $C(1)$) corresponds to the trivial
(resp. unique irreducible nontrivial) bundle on the regular nilpotent orbit.
$G-$equivariant bundles on the point bijectively correspond to representations
of $G$. Under this identification the element $C(n), n\ge 2$ correspond to
irreducible $SL_2-$representation with highest weight $n-2$. We see that Conjectures
1---4 hold in this case. Moreover it is easy to check that all tilting
Conjectures are true in this case.

\subsection{Regular nilpotent orbit} The support of an element $C(\lambda)$ is the
full nilpotent cone if and only if $\lambda$ is a {\em minimal} weight. In this case
$C(\lambda)=AJ(\lambda)=(-v)^{\nu -\nu_{\lambda}}e^{\lambda}$. This fits nicely
with results of Graham \cite{Gr} (see also \cite{Gi}) who computed the $G-$module
structure of the ring of functions on universal cover $\tilde \cO$ of the regular
nilpotent orbit. Namely, Graham proved the following equality in $K_{\cG}(\cN)$:
$$
[\BC [\tilde \cO]]=\sum_{\lambda \; \mbox{is minimal}}v^{\nu -\nu_{\lambda}}e^{\lambda}.
$$

\subsection{Lowest cell} The support of an element $C(\lambda)$ should be a 
point
if and only if $\lambda -2\varrho$ is dominant. In this case $C(\lambda)$ should
correspond to the irreducible representation of $G$ with highest weight 
$\lambda -2\varrho$.
Using the Koszul complex we see that this is equivalent to the equality
$$
C(\lambda)=v^{-\nu}e^{\lambda -2\varrho}\prod_{\alpha \in R_+}(e^0-v^2e^{\alpha}) \eqno (*)
$$
where $R_+$ is the set of positive roots. This formula should be understood
as follows: first we make all multiplications and then interpret the result as
an element of $K_{\cG}(\cN)$. The reader should be aware that the map $s_*: K_{\cG}(T^*\cB)\to
K_{\cG}(\cN)$ {\em is not} multiplicative (and moreover the group $K_{\cG}(\cN)$
has no natural multiplicative structure). We say that $\lambda \in X_+$ is {\em
very dominant} if for any subset $J\subset R_+$ the weight $\lambda +
\sum_{\alpha \in J}\alpha$ is dominant. One can show that the right hand side of (*)
is a selfdual element of $K_{\cG}(\cN)$. Now it is clear from the definitions that
formula (*) is true for any very dominant $\lambda$. It would be interesting to
prove the formula (*) in general, the most interesting case being $\lambda =2\varrho$.
I checked this formula for groups of rank 2.

\subsection{McGovern formula} Conjecture 4 is very easy to check in each
particular case thanks to McGovern's formula \cite{Mc1} for $G-$structure of the
ring of functions on nilpotent orbits. We restate this formula as follows.
The Dynkin diagram of a nilpotent orbit determines a grading of the set of positive
roots (this grading is additive and the gradation of a simple root is its label
in the Dynkin diagram). Let $R_{+0}\subset R_+$ (resp. $R_{+1}\subset R_+$) be the
set of positive roots with gradation 0 (resp. 1). The McGovern formula is
the specialized at $v=1$ right hand side of the following version
of Conjecture 4:
$$
C(\lambda_{\cO})=v^{-a(\cO)}\prod_{\alpha \in R_{+0}\cup R_{+1}}(e^0-v^2e^{\alpha}).
$$
Here $a(\cO)=|R_{+0}|+\frac{1}{2}|R_{+1}|$. We remark that the right hand side of
this formula is selfdual (as mentioned after the Conjecture 4). We checked
that this formula works for groups of rank 2, for the subregular nilpotent 
orbit, 
and in some other cases. As a consequence of this formula we obtain a 
combinatorial formula for $\lambda_{\cO}$ and explicit formula for dominant 
distinguished involutions in $\hat W^a$. It would be very interesting to find 
such formulas for other cases mentioned after the Conjecture 4.

In the special case of the group $SL_n$ the formula for $\lambda_{\cO}$ can be 
described (following a remark by David Vogan) quite explicitly as follows. 
Let the sizes of Jordan blocks of an element $e\in \cO$ be given by partition 
$p=p_1\ge p_2\ge \ldots$. Let $p'=p'_1\ge p'_2\ge \ldots$ be the dual partition
and let $\cO'$ be the nilpotent orbit consisting of matrices with Jordan blocks
of sizes $p'_1, p'_2, \ldots$. Let us consider the Dynkin diagram of $\cO'$ as
a weight $\lambda$ for $SL_n$. Then $\lambda_{\cO}=\lambda$. It would be 
interesting to find similar combinatorial formulas for other classical groups.

\subsection{$SL_3$}
\label{sl3}
Let $\omega_1$ and $\omega_2$ denote fundamental weights of $SL_3$. The weights
of $SL_3$ not covered by previous discussion are weights "near the walls"
$n\omega_i, n\ge 2, i=1,2$ and $\omega_1+\omega_2+n\omega_i, n\ge 1, i=1,2$
corresponding to the subregular nilpotent orbit of $SL_3$. One checks that in
this case all our Conjectures are true. Here we only consider the weight
$\lambda =3\omega_1$. One calculates $C(3\omega_1)=AJ(3\omega_1)+(v^{-1}+v^{-3})
AJ(\omega_1+\omega_2)+v^{-2}AJ(0)=v^2e^{3\omega_1}-(1+v^2)e^{\omega_1+\omega_2}+
v^{-2}e^0$. Further, $\Gamma (C(3\omega_1))=v^{-2}V_0-v^2V_{3\omega_2}-v^4\ldots$
(cf. with last paragraph of \cite{O}).
In particular, $C(3\omega_1)$ is not of the form $[\cS]$ for some sheaf $\cS$.

\subsection{Subregular nilpotent orbits} Let $\cO$ be the subregular nilpotent
orbit. For any simple root $\alpha_i$ let $P_i$ be the corresponding parabolic
subgroup. As it is well known for any {\em short} simple root $\alpha_i$ the 
moment map $T^*(G/P_i)\to \fg$ is a resolution of singularities of the closure 
$\bar \cO$, see e.g. \cite{Br}. We get that $\lambda_{\cO}$ is the unique short
dominant root and
$C(\lambda_{\cO})=v^{-1}e^0-ve^{\alpha_i}$.

Let $G_{ad}$ be the adjoint group of the same type as $G$. A nontrivial
$G_{ad}-$equivariant irreducible bundle on $\cO$ exists if and only if
$G$ is not simply laced. In cases of types $B_n, C_n (n\ge 2), F_4$ such a
bundle $\cV$ is unique (the case of type $G_2$ is considered below).  
For any {\em long} simple
root $\alpha_j$ the image of the moment map $T^*(G/P_j)$ is $\bar \cO$ and this
map is generically two to one. We deduce that the weight $\lambda_{\cV}$
corresponding to the bundle $\cV$ is the unique long dominant root and 
$C(\lambda_{\cV})=v(e^{\alpha_i}-e^{\alpha_j})$, where $\alpha_i$ is a
short simple root.


\subsection{The subregular nilpotent orbit in $G_2$} The fundamental group of
the subregular nilpotent orbit $\cO$ for the group $G$ of type $G_2$ is the 
symmetric group in
three letters $S_3$. It is not hard to guess what weights should correspond
to irreducible representations of $S_3$. Let $\omega_1$ and $\omega_2$ be the
fundamental weights for $G_2$ such that $\dim V_{\omega_1}=14$ and
$\dim V_{\omega_2}=7$. The weight $\omega_2$ (resp. $\omega_1,\; 2\omega_2$)
corresponds to the trivial representation of $S_3$ (resp. the irreducible
two dimensional, the sign representation). We have
$$
C(\omega_2)=AJ(\omega_2)+v^{-1}AJ(0)=v^{-1}(e^0-v^6e^{\omega_2}),
$$
$$
C(\omega_1)=AJ(\omega_1)-v^{-4}AJ(\omega_2)=v^{-1}(v^2e^{\omega_2}-v^6e^{\omega_1}),
$$
$$
C(2\omega_2)=AJ(2\omega_2)-v^{-2}AJ(\omega_1)=v^{-1}(v^4e^{\omega_1}-v^6e^{2\omega_2}).
$$
This implies the following formula for the image of the trivial bundle $\BC [\tilde \cO ]$
on the universal cover $\tilde \cO$ in the group $K_{\cG}(\cN)$:
$$
[\BC [\tilde \cO ]]=C(\omega_2)+2C(\omega_1)+C(2\omega_2)=v^{-1}(e^0+2v^2e^{\omega_2}
+v^4e^{\omega_1}-v^6e^{\omega_2}-2v^6e^{\omega_1}-v^6e^{2\omega_2}).
$$
In particular we obtain a (conjectural) formula for graded multiplicities of simple $G-$modules
in the function ring of $\tilde \cO$. Fortunately, another formula for these
multiplicities is available in the literature thanks to the work of
McGovern \cite{Mc2}. We checked that McGovern's formula and ours
are equivalent.

\end{document}